\documentstyle[11pt]{article}
\newcommand{\ba}{\noindent $\begin{array}}
\newcommand{\ea}{\end{array}$}
\newcommand{\be}{\begin{equation}}
\newcommand{\ee}{\end{equation}}
\newcommand{\bd}{\begin{displaymath}}
\newcommand{\ed}{\end{displaymath}}
\newcommand{\beq}{\begin{eqnarray*}}
\newcommand{\eeq}{\end{eqnarray*}}
\newcommand{\beqn}{\begin{eqnarray}}
\newcommand{\eeqn}{\end{eqnarray}}

\begin{document}

\pagestyle{plain}

\begin{center}

{\large \bf COMPUTATIONAL EXPERIMENTS WITH  ABS ALGORITHMS
\vskip2mm FOR OVERDETERMINED LINEAR SYSTEMS}
\footnote{Work  supported by grant GA CR 201/00/0080
and by MURST 1997 and 1999 
 Cofinanziamento}

\end{center}

\begin{center}
{\bf E. Bodon
\footnote{Department of Mathematics, University of Bergamo, Bergamo 24129,
Italy (bodon@unibg.it)},
A. Del Popolo
\footnote{Department of Mathematics, University of Bergamo, Bergamo 24129,
Italy (delpopolo@unibg.it)},
L. Luk\v{s}an
\footnote{Institute of Computer Science, Academy of Sciences of the Czech
Republic,
Pod vod\'arenskou v\v e\v z\'\i\ 2, 182 07 Prague 8, Czech Republic
(luksan@uivt.cas.cz)} and
E. Spedicato
\footnote{Department of Mathematics, University of Bergamo, Bergamo 24129,
Italy (emilio@unibg.it)}}
\end{center}

\vspace{2mm}

\section{Introduction}

In this report, we present numerical experiments with two particular 
ABS algorithms:

\begin{itemize}

\item[(1)] The Huang or modified Huang algorithm,

\item[(2)] The implicit QR algorithm,

\end{itemize}

\noindent These methods are used for solving the following problem:

\begin{itemize}

\item[(a)] Finding the least-squares solution of an overdetermined linear
system $A x = b$, where $A \in R^{m,n}$, $b \in R^m$, $x \in R^n$ and 
$m \geq n$.

\end{itemize}

\noindent The considered algorithms belong to the 
scaled  ABS class of algorithms for solving linear systems,
 that is defined by the following scheme:

\begin{itemize}

\item[(A)] Let $x_1 \in R^n$ be arbitrary and $H_1 \in R^{n,n}$ be  
nonsingular arbitrary. Set $i = 1$.

\item[(B)] Compute the residual $r_i = A x_i - b$. If $r_i = 0$, then 
stop, $x_i$ solves the problem. Otherwise compute $s_i = H_i A^T v_i$, 
where $v_i \in R^n$ is arbitrary, save that $v_1, \dots, v_i$ are linearly 
independent. If $s_i \neq 0$, then  go  to (C). If $s_i = 0$ and 
$r_i^T v_i = 0$, then set $x_{i+1} = x_i$, $H_{i+1} = H_i$ and go to (F). 
If $s_i = 0$ and $r_i^T v_i \neq 0$, then stop, the system is incompatible.

\item[(C)] Compute the search vector $p_i$ by
$$
p_i = H_i^T z_i,
$$
where $z_i \in R^n$ is arbitrary, save that $z_i^T H_i A^T v_i \neq 0$.

\item[(D)] Update the estimate of the solution by
$$
x_{i+1} = x_i - \alpha_i p_i,
$$
where the stepsize $\alpha_i$ is given by
$$
\alpha_i = r_i^T v_i / p_i^T A^T v_i. 
$$

\item[(E)] Update the Abaffian matrix by
$$
H_{i+1} = H_i - H_i A^T v_i w_i^T H_i / w_i^T H_i A^T v_i,
$$
where $w_i \in R^n$ is arbitrary, save that $w_i^T H_i A^T v_i \neq 0$.

\item[(F)] If $i = m$, then stop, $x_{i+1}$ solves the problem. Otherwise 
increment the index $i$ by one and go to (B).

\end{itemize}

From (E), it follows by induction that $H_{i+1} A^T V_i = 0$ and  
$H_{i+1}^T W_i = 0$, where $V_i = [v_1, \dots, v_i]$ and 
$W_i = [w_1, \dots, w_i]$.
One can show that the implicit factorization $V_i^T AP_i = L_i$ holds, where
$P_i = [p_1, \dots, p_i]$ and $L_i$ is nonsingular lower triangular.
Moreover the general solution of the scaled subsystem $V_i^T A x =V_i^T b$ 
can be expressed in the form
\be
x = x_i + H_i^T q, 
\label{0}
\ee
where $q \in R^n$ is arbitrary (see \cite{ab2} for the proof).

The basic ABS class is the subclass of the scaled ABS class
obtained by taking $v_i = e_i$, $e_i$ being the $i$-th unitary vector 
($i$-th column of the unit matrix). In this case, residual $r_i = A x_i - b$ 
need not be computed in (B), $i$-th element $r_i^T e_i = a_i^T x_i - b_i$ 
suffices. 

This report is organized as follows. In Section 2, a short description of 
individual algorithms is given. Section 3 contains some details concerning
test matrices and numerical experiments and diuscusses the
numerical results. The Appendix contains listing of all numerical results.
For other numerical results on ABS methods see [1,11]. For a listing
of the used ABS codes see [12].
\section{The tested ABS  algorithms for overdetermined linear systems}

In this report, we deal with two basic algorithms, one belonging to the
basic ABS class and one belonging to the scaled ABS class. To simplify
description, we will assume that $A \in R^{m,n}$, 
$m \leq n$, has full row rank so that $s_i \neq 0$ in (B).

The Huang algorithm is obtained by the parameter choices $H_1 = I$, 
$v_i = e_i$, $z_i = a_i$, $w_i = a_i$. Therefore
\be
p_i = H_i a_i
\label{1}
\ee
and
\be
H_{i+1} = H_i - p_i p_i^T / a_i^T p_i,
\label{2}
\ee
From (\ref{1}) and (\ref{2}), it follows by induction that 
$p_i \in Range(A_i)$ and 
\be
H_{i+1} = I - P_i D_i^{-1} P_i^T. 
\label{3}
\ee
where $A_i = [a_1, \dots, a_i]$, $P_i = [p_1, \dots, p_i]$ and 
$D_i = diag(a_1^Tp_1, \dots, a_i^Tp_i)$. Moreover $H_i$ is symmetric,
positive semidefinite and idempotent (it is the orthogonal projection
matrix into $Null(A_{i-1})$). Since the requirement 
$p_i \in Null(A_{i-1})$ is crucial, we can improve orthogonality by 
iterative refinement $p_i = H_i^j a_i$, $j > 1$ (usually $j=2$), 
obtaining the modified Huang algorithm.  

The Huang algorithm can be used for finding the minimum-norm solution to 
the compatible underdetermined system $A x = b$, i.e. for minimizing 
$\| x \|$ s.t. $A x = b$. To see this, we use the Lagrangian function 
and convexity  of $\| x \|$. Then $x$ is a required solution if and only 
if $x = A^T u$  for some $u \in R^m$ or $x \in Range(A^T)$. But  
$$
x_{m+1} = x_1 - \sum_{i=1}^m \alpha_i p_i 
$$ 
by (D) and $p_i \in Range(A_i) \subset Range(A^T)$, so that if $x_1 = 0$, 
then $x_{m+1} \in Range(A^T)$. 

A short description of two versions of the Huang and modified Huang
algorithms follows.

\vspace{4mm}

\noindent {\bf Algorithm 1}

\noindent (Huang and modified Huang, formula (\ref{2})).

\noindent Set $x_1 = 0$ and $H_1 = I$.

\noindent {\bf For} $i=1$ {\bf to} $m$ {\bf do}

\noindent Set $p_i = H_i a_i$ (Huang) or 

\noindent $p_i = H_i (H_i a_i)$ (modified Huang),

\noindent $d_i = a_i^T p_i$ and $x_{i+1} = x_i - ((a_i^T x_i - b_i)/d_i) p_i$.

\noindent If $i < m$, then set $H_{i+1} = H_i - p_i p_i^T / d_i$.

\noindent {\bf end do}

\vspace{4mm}

\noindent {\bf Algorithm 2}

\noindent (Huang and modified Huang, formula (\ref{3})).

\noindent Set $x_1 = 0$ and $P_0$ empty.

\noindent {\bf For} $i=1$ {\bf to} $m$ {\bf do}

\noindent Set $p_i = (I - P_{i-1} D_{i-1}^{-1} P_{i-1}^T) a_i$ (Huang) or 

\noindent $p_i = (I - P_{i-1} D_{i-1}^{-1} P_{i-1}^T) 
		 (I - P_{i-1} D_{i-1}^{-1} P_{i-1}^T) a_i$  (modified Huang),
		 
\noindent $d_i = a_i^T p_i$ and $x_{i+1} = x_i - ((a_i^T x_i - b_i)/d_i) p_i$.

\noindent If $i < m$, then set $P_i = [ P_{i-1}, p_i ]$. 

\noindent {\bf end do}

\vspace{4mm}

The implicit QR algorithm is obtained by the parameter choices $H_1 = I$, 
$v_i = A p_i$, $z_i = e_i$, $w_i = e_i$.  Since $z_i = e_i$ and $w_i = e_i$,
the matrices $H_i$ and $P_i$ have the same structure as that in the implicit 
LU algorithm. Using the implicit factorization property $V_i^T AP_i = L_i$,
we can see that $V_{i-1}^T v_i = V_{i-1}^T Ap_i = 0$ so that the vectors 
$v_j$, $j = 1, \dots, i$, are mutually orthogonal.  

The implicit QR algorithm can be used for finding the least-squares solution
of an overdetermined linear system $A x = b$, where $A \in R^{m,n}$, 
$b \in R^m$, $x \in R^n$ and $m \geq n$, which is obtained after at most 
$n$  steps. To see this, we recall that algorithms from the scaled ABS class 
solve the scaled system $V_n^T A x = V_n^T b$. Since $V_n = A P_n$, where $P_n$
is square nonsingular, the condition $P_n^T A^T A x = P_n^T A^T b$ implies 
$A^T A x = A^T b$, which defines the least-squares solution of the 
overdetermined system $A x = b$.

A short description of the implicit QR algorithm follows.

\vspace{4mm}

\noindent {\bf Algorithm 3}

\noindent (Implicit QR).

\noindent Set $x_1 = 0$, $r_1 = -b$ and $H_1 = I$.

\noindent {\bf For} $i=1$ {\bf to} $m$ {\bf do}

\noindent Set $p_i = H_i^T e_i$, $v_i = A p_i$  

\noindent (only $(i-1)(n-i+1)$ nonzero elements of $H_i$ is used),

\noindent $\alpha_i =r_i^T v_i / v_i^T v_i$, $x_{i+1} = x_i - \alpha_i p_i$ 
	  and $r_{i+1} = r_i - \alpha_i v_i$.

\noindent (only $i$ nonzero elements of $x_i$ are updated).

\noindent If $i < m$, then set $s_i = H_i^T A^T v_i$ and 
	  $H_{i+1} = H_i - s_i p_i^T / v_i^T v_i$

\noindent (only $i(n-i)$ nonzero elements of $H_{i+1}$ are updated).

\noindent {\bf end do}

\vspace{4mm}

We now consider the problem of solving overdetermined linear systems,
in the least squares sense.
Consider the linear system $A x = b$, where $A \in R^{m,n}$, $b \in R^m$, 
$x \in R^n$ and $m \geq n$. This system can be solved in the least
squares sense
in $n$ steps
by the inplicit QR algorithm as was shown above. Another possibility is
based on the application of the Huang algorithm. Since the least-squares
solution has to satisfy the normal equation 
$$ 
A^T A x = A^T b,
$$
it can be obtained from the solution of the following augmented 
system .
\be
A x = y, 
\label{5}
\ee
\be
A^T y = A^T b. 
\label{6}
\ee
Since $y \in Range(A)$ by (\ref{5}), it is a minimum-norm solution of
(\ref{6}) and it can be obtained by the Huang algorithm. Having $y$,
the compatible system (\ref{5}) can be solved by any ABS algorithm. 
Moreover, using the implicit factorization property $A^T P_n = L_n$,
we can write
$$
L_n^T x = P_n^T A x = P_n^T y = \tilde{b} 
$$
which as a system with a lower triangular matrix can be easily 
solved by the back substitution. The lower triangular matrix $L_n$ can
be obtained as a by-product of the Huang algorithm. If we use 
formula (\ref{3}), then $d_i$ is $i$-th diagonal element of $L_n$ and
$P_{i-1}^T \tilde{a}_i$ contains the other $i$-th-column elements of $L_n$
($\tilde{a}_i$ is $i$-th column of the matrix $A$).

Matrix $L_n$ has not to be stored since it can be reconstructed from
columns of the matrix $P_n$. However, the back substitution has 
to be realized in a slightly different way in this case.     

A short description of two versions of the Huang algorithm for least-squares 
solution of overdetermined systems follows.

\vspace{4mm}

\noindent {\bf Algorithm 4}

\noindent (Huang and modified Huang for least-squares with stored $L_n$).

\noindent Set $x_1 = 0$ and $P_0$ empty.

\noindent {\bf For} $i=1$ {\bf to} $n$ {\bf do}

\noindent If $i > 1$, then compute $g_i = P_{i-1}^T \tilde{a}_i$ and 
	  copy it into the lower triangular matrix $L_n$.

\noindent Set $p_i = \tilde{a}_i - P_{i-1} D_{i-1}^{-1} g_i$ (Huang) or 

\noindent $p_i = (I - P_{i-1} D_{i-1}^{-1} P_{i-1}^T) 
		 (\tilde{a}_i - P_{i-1} D_{i-1}^{-1} g_i)$  (modified Huang).
		 
\noindent Set $d_i = \tilde{a}_i^T p_i$ and copy it into the lower 
	  triangular matrix $L_n$. Set $\tilde{b}_i = b^T p_i$.
	  
\noindent If $i < n$, then set $P_i = [ P_{i-1}, p_i ]$. 

\noindent {\bf end do}

\noindent Solve the triangular system $L_n^T x = \tilde{b}$.

\vspace{4mm}

\noindent {\bf Algorithm 5}

\noindent (Huang and modified Huang for least-squares without stored $L_n$).

\noindent Set $x_1 = 0$ and $P_0$ empty.

\noindent {\bf For} $i=1$ {\bf to} $n$ {\bf do}

\noindent Set $p_i = (I - P_{i-1} D_{i-1}^{-1} P_{i-1}^T) \tilde{a}_i$ (Huang) or 

\noindent $p_i = (I - P_{i-1} D_{i-1}^{-1} P_{i-1}^T) 
		 (I - P_{i-1} D_{i-1}^{-1} P_{i-1}^T) \tilde{a}_i$  (modified Huang).
		 
\noindent Set $d_i = \tilde{a}_i^T p_i$. 
	  
\noindent If $i < n$, then set $P_i = [ P_{i-1}, p_i ]$. 

\noindent {\bf end do}

\noindent Set $f_n = b$. 

\noindent {\bf For} $i=n$ {\bf to} $1$ {\bf do}

\noindent Set $x_i = p_i^T f_i / d_i$.

\noindent If $i > 1$, then set $f_{i-1} = f_i - x_i \tilde{a}_i$  

\noindent {\bf end do}

\vspace{4mm}

\section{The computational experiments}

Performance of ABS algorithms has been tested by using several types of 
ill-conditioned matrices. These matrices can be classified in the following
way. The first letter distinguishes matrices with integer 'I' and real 
'R' elements, both actually stored as reals in double precision arithmetic.      
The second letter denotes randomly generated matrices 'R' or matrices 
determined by an explicit formula 'D'. For randomly generated matrices,
a number specifying the interval for the random number generator follows,
while the name of matrices determined by the explicit formula contains
the formula number (F1, F2, F3). The last letter of the name denotes
a way for obtaining ill-conditioned matrices: 'R' - matrices with nearly
dependent rows, 'C' - matrices with nearly dependent columns, 'S' - nearly
singular symmetric matrices, 'B' -  both matrices in KKT system ill-conditioned.
More specifically:

\vspace{0.2cm}

\noindent \parbox{1.5cm}{IR500 } \parbox[t]{14.4cm}{Randomly generated matrices 
with integer elements uniformly distributed in the interval [-500,500].}

\noindent \parbox{1.5cm}{IR500R} \parbox[t]{14.4cm}{Randomly generated matrices 
with integer elements uniformly distributed in the interval [-500,500] perturbed 
in addition to have two rows nearly dependent.}

\noindent \parbox{1.5cm}{IR500C} \parbox[t]{14.4cm}{Randomly generated matrices 
with integer elements uniformly distributed in the interval [-500,500] perturbed 
in addition to have two columns nearly dependent.}

\noindent \parbox{1.5cm}{RR100 } \parbox[t]{14.4cm}{Randomly generated matrices 
with real elements uniformly distributed in the interval [-100,100]}.

\noindent \parbox{1.5cm}{IDF1  } \parbox[t]{14.4cm}{Matrices with elements 
$a_{ij} = |i-j|$, $1 \leq i \leq m$, $1 \leq j \leq n$ (Micchelli-Fiedler matrix).} 

\noindent \parbox{1.5cm}{IDF2  } \parbox[t]{14.4cm}{Matrices with elements 
$a_{ij} = |i-j|^2$, $1 \leq i \leq m$, $1 \leq j \leq n$.}

\noindent \parbox{1.5cm}{IDF3  } \parbox[t]{14.4cm}{Matrices with elements 
$a_{ij} = |i+j-(m+n)/2|$, $1, \leq i \leq m$, $1 \leq j \leq n$.}

\noindent \parbox{1.5cm}{IR50  } \parbox[t]{14.4cm}{Randomly generated matrices
with integer elements uniformly distributed in the interval [-50,50].}

\vspace{0.2cm}

Matrices with linearly dependent rows were obtained in the following way. 
The input data contain four integers which specify two row indices $i_1$, 
$i_2$, one column index $i_3$ and one exponent $i_4$. Then the row $a_{i_1}$
is copied into $a_{i_2}$. Furthermore $a_{i_1 i_3}$ is set to zero and
$a_{i_2 i_3}$ to $2^{-i_4}$. Similar procedures are used for columns 
and symmetric matrices.

Solution vectors were  generated randomly
with integer or real elements uniformly distributed in the interval [-10,10]. 
Right hand sides for overdetermined systems
were obtained by the following way. First, vector $\tilde{b}$ was 
generated randomly with integer or real elements uniformly distributed in the 
interval [-10,10]. Then its first element together with elements in the first 
row of the matrix $A$ were redefined by the formulas $\tilde{b}_1 = -1$ and 
$a_{1j} = \sum_{i=2}^n a_{ij} \tilde{b}_j$ so that $A^T \tilde{b} = 0$. The 
right-hand side vector was determined by the formula $b = \tilde{b} + A x^{\star}$.
Since $A^T A x^{\star} = A^T (b - \tilde{b}) = A^T b$, the normal equation
is satisfied and $x^{\star}$ is a least squares solution of the system $A x = b$.
Vector $\tilde{b}$ was generated nonzero, so that the system $A x = b$ is 
incompatible. 

We have tested the following particular algorithms:

\vspace{0.2cm}

\noindent \parbox{2.5cm}{impl.qr5  } \parbox[t]{13.4cm}{The implicit QR algorithm: 
Algorithm 5 (subroutine alg5d.f).}

\noindent \parbox{2.5cm}{huang6    } \parbox[t]{13.4cm}{The Huang algorithm for 
overdetermined systems: Algorithm 6 (subroutine alg6d.f).}.

\noindent \parbox{2.5cm}{mod.huang6} \parbox[t]{13.4cm}{The modified Huang algorithm  
for overdetermined systems: Algorithm 6 (subroutine alg6d.f).}

\noindent \parbox{2.5cm}{huang7    } \parbox[t]{13.4cm}{The Huang algorithm for 
overdetermined systems: Algorithm 7 (subroutine alg7d.f).}

\noindent \parbox{2.5cm}{mod.huang7} \parbox[t]{13.4cm}{The modified Huang algorithm  
for overdetermined systems: Algorithm 7 (subroutine alg7d.f).}

\noindent \parbox{2.5cm}{qr lapack } \parbox[t]{13.4cm}{Subroutine DGELS  
from the LAPACK package.}

\noindent \parbox{2.5cm}{svd lapack} \parbox[t]{13.4cm}{Subroutine DGELSS 
from the LAPACK package.}

\noindent \parbox{2.5cm}{gqr lapack} \parbox[t]{13.4cm}{Subroutine DGELSX 
from the LAPACK package.}

\vspace{0.2cm}

\noindent Notice that ABS algorithms were implemented in their basic form without
partitioning into block or other special adjustements serving for speed increase 
as done in the LAPACK software.

Detailed results of computational experiments are presented in the Appendix. For each selected
problem, the type of matrix and the dimension is given. Furthermore, both the solution 
and the residual errors together with the detected rank and the computational time 
are given for each tested algorithm. Computational experiments were performed on a 
Digital Unix Workstation in the double precision arithmetic (machine epsilon equal
to about $10^{-16}$).

We have tested overdetermined systems with $n >> m/2$, $n = m/2$ and 
$n << m/2$. These systems were solved by using Algorithms 5 - Algorithm 7 together with 
explicit QR and SVD decomposition based methods taken  from the LAPACK package.

 The following  tables give
synthetic results for the 21 tested problems, the number at the intersection
of the $i$-th row with the $k$-th column indicating how many times the
algorithm at the heading of the $i$-th row gave a lower error 
than the algorithm at the heading of the $k$-th row (in case there is
a second number, this counts the number of cases when difference  was
less that one percent).
\vskip3mm
\begin{verbatim}

   solution error  -  21  overdetermined linear systems

   methods       huang  mod.  huang  mod.   impl.  qr   svd   gqr
   "             6     huang  7     huang  qr5   lap   lap   lap 
   "                     6             7          ack   ack   ack        total

   huang6               4      0/21  4      5     9     2     2         26/21
   mod.huang6   17            17     1/12  13/6  16    11     8         83/18
   huang7        0/21   4            4      5     9     2     2         26/21
   mod.huang7   17      8/12  17           14/6  18    11    10         95/18
   impl.qr5     16      2/6   16     1/6         10     8     6         59/12
   qr lapack    12      5     12     3     11           6/4   9         58/4
   svd lapack   19     10     19    10     13    11/4         4/12      86/16
   gqr lapack   19     13     19    11     15    12     5/12            94/12
\end{verbatim}
\newpage
\begin{verbatim}


   residual error  -  21  overdetermined linear systems

   methods       huang  mod.   huang  mod.  impl.  qr    svd   gqr
   "             6     huang   7     huang  qr5   lap   lap   lap
   "                    6             7           ack   ack   ack      total

   huang6              11      4     7     17     16     9    13        77
   mod.huang6    10            6/2   9     15/1   13    10    11/1      74/4
   huang7        17    13/2          9     16     16    12    12        95/2
   mod.huang7    14    12     12           17     14/1  13    13        95/1
   impl.qr5       4     5/1    5     4             8     7     7        40/1
   qr lapack      5     8      5     6/1   13            9     9        55/1
   svd lapack    12    11      9     8     14     12          13        79
   gqr lapack     8     9/1    9     8     14     12     8              68/1

\end{verbatim}

From the results in the Appendix and the above tables we can state the following
conclusions:

\begin{itemize}

\item[(1)] When solving well-conditioned overdetermined systems, the explicit QR algorithms 
based on the Householder orthogonalization are usually faster than ABS methods tested 
(especially if $n >> m/2$). Therefore, we cannot recommend the later for solving standard 
problems. 

\item[(2)] Interesting results were obtained for extremely ill-conditioned problems. 
The modified Huang and the implicit QR algorithms failed to solve systems with the 
matrix IDF2. Such systems was not solved by simple explicit QR methods as well. On the 
other hand, the modified Huang and the implicit QR algorithms found solutions of systems 
with the matrix IDF3 extremely fast and, moreover, they were able to determine the numerical 
rank correctly. Performance of these algorithms in that particular case is remarkable, 
they are at least 20 times faster then the best LAPACK routines. 

\item[(3)] The LAPACK routine DGELSS based on the SVD decomposition is extremely slow 
(especially if $n >> m/2$), not suitable for solving our problems (it can be substituted by
the much faster routine DGELSX). 

\item[(4)] In term of solution error mod-huang7 and gqr are the best; in term
of residual error huang7 and mod.huang7 are the best.

\end{itemize}

\newpage

\section{Appendix: Test results on overdetermined linear systems}

\small

\begin{verbatim}           

 matrix  dimension    method       solution  residual   rank     time
	   m    n                    error     error
 --------------------------------------------------------------------

 IR500   1050  950    huang6       0.17D-01  0.23D-14    950    32.00
 IR500   1050  950    mod.huang6   0.78D-10  0.20D-15    950    63.00
 IR500   1050  950    huang7       0.17D-01  0.16D-14    950    30.00
 IR500   1050  950    mod.huang7   0.79D-10  0.85D-15    950    60.00
 IR500   1050  950    impl.qr5     0.51D-09  0.12D-13    950    57.00
 IR500   1050  950    qr lapack    0.84D-09  0.83D-14    950    16.00
 IR500   1050  950    svd lapack   0.31D-08  0.80D-14    950   119.00
 IR500   1050  950    gqr lapack   0.28D-08  0.28D-14    950    26.00
condition number:     0.42D+09
 
 IR500   1400  700    huang6       0.13D-02  0.51D-14    700    23.00
 IR500   1400  700    mod.huang6   0.41D-11  0.12D-13    700    49.00
 IR500   1400  700    huang7       0.13D-02  0.82D-14    700    24.00
 IR500   1400  700    mod.huang7   0.39D-11  0.15D-14    700    47.00
 IR500   1400  700    impl.qr5     0.33D-10  0.31D-13    700    40.00
 IR500   1400  700    qr lapack    0.54D-10  0.56D-14    700    15.00
 IR500   1400  700    svd lapack   0.54D-10  0.18D-13    700    67.00
 IR500   1400  700    gqr lapack   0.23D-09  0.77D-14    700    23.00
condition number:     0.65D+08
\end{verbatim}           

\newpage

Test results for overdetermined linear systems - continued

\begin{verbatim}           
 matrix  dimension    method       solution  residual   rank     time
	   m    n                    error     error
 --------------------------------------------------------------------

 IR500   2000  400    huang6       0.35D-03  0.95D-15    400    12.00
 IR500   2000  400    mod.huang6   0.96D-12  0.78D-15    400    21.00
 IR500   2000  400    huang7       0.35D-03  0.41D-15    400    11.00
 IR500   2000  400    mod.huang7   0.87D-12  0.33D-15    400    21.00
 IR500   2000  400    impl.qr5     0.55D-11  0.11D-13    400    18.00
 IR500   2000  400    qr lapack    0.18D-10  0.46D-13    400     7.00
 IR500   2000  400    svd lapack   0.18D-10  0.54D-13    400    20.00
 IR500   2000  400    gqr lapack   0.65D-10  0.10D-12    400    12.00
 condition number:    0.24D+08

 IR500C  1050  950    huang6       0.23D-01  0.15D-14    950    33.00
 IR500C  1050  950    mod.huang6   0.11D-01  0.19D-14    950    63.00
 IR500C  1050  950    huang7       0.23D-01  0.36D-16    950    30.00
 IR500C  1050  950    mod.huang7   0.11D-01  0.11D-14    950    60.00
 IR500C  1050  950    impl.qr5     0.12D-01  0.40D-12    950    57.00
 IR500C  1050  950    qr lapack    0.45D-02  0.36D-14    950    17.00
 IR500C  1050  950    svd lapack   0.69D+00  0.41D-14    949   114.00
 IR500C  1050  950    gqr lapack   0.69D+00  0.19D-14    949    26.00
condition number:     0.25D+14

 IR500C  1400  700    huang6       0.27D+00  0.11D-14    700    24.00
 IR500C  1400  700    mod.huang6   0.21D-05  0.34D-15    700    49.00
 IR500C  1400  700    huang7       0.27D+00  0.62D-15    700    24.00
 IR500C  1400  700    mod.huang7   0.21D-05  0.23D-15    700    47.00
 IR500C  1400  700    impl.qr5     0.12D+01  0.65D-14    700    29.00
 IR500C  1400  700    qr lapack    0.35D-06  0.11D-15    700    15.00
 IR500C  1400  700    svd lapack   0.35D-06  0.80D-15    700    66.00
 IR500C  1400  700    gqr lapack   0.11D-05  0.32D-14    700    22.00
condition number:     0.51D+12
 
 IR500C  2000  400    huang6       0.11D-01  0.31D-14    400    11.00
 IR500C  2000  400    mod.huang6   0.33D-03  0.18D-14    400    21.00
 IR500C  2000  400    huang7       0.11D-01  0.50D-15    400    11.00
 IR500C  2000  400    mod.huang7   0.33D-03  0.37D-15    400    21.00
 IR500C  2000  400    impl.qr5     0.80D-03  0.23D-13    400    18.00
 IR500C  2000  400    qr lapack    0.10D-02  0.37D-15    400     7.00
 IR500C  2000  400    svd lapack   0.70D+00  0.55D-15    399    19.00
 IR500C  2000  400    gqr lapack   0.70D+00  0.19D-13    399    11.00
condition number:     0.52D+13
\end{verbatim}           

\newpage

Test results for overdetermined linear systems - continued

\begin{verbatim}           
 matrix  dimension    method       solution  residual   rank     time
	   m    n                    error     error
 --------------------------------------------------------------------

 RR100   1050  950    huang6       0.59D-10  0.21D-14    950    32.00
 RR100   1050  950    mod.huang6   0.96D-14  0.25D-14    950    63.00
 RR100   1050  950    huang7       0.59D-10  0.93D-15    950    30.00
 RR100   1050  950    mod.huang7   0.85D-14  0.58D-15    950    61.00
 RR100   1050  950    impl.qr5     0.82D-14  0.31D-14    950    57.00
 RR100   1050  950    qr lapack    0.79D-14  0.19D-15    950    16.00
 RR100   1050  950    svd lapack   0.20D-13  0.30D-14    950   158.00
 RR100   1050  950    gqr lapack   0.80D-14  0.26D-14    950    26.00
condition number:     0.35D+04

 RR100   1400  700    huang6       0.25D-11  0.47D-12    700    24.00
 RR100   1400  700    mod.huang6   0.44D-14  0.99D-12    700    48.00
 RR100   1400  700    huang7       0.25D-11  0.99D-12    700    24.00
 RR100   1400  700    mod.huang7   0.24D-14  0.85D-13    700    47.00
 RR100   1400  700    impl.qr5     0.35D-14  0.54D-12    700    40.00
 RR100   1400  700    qr lapack    0.28D-14  0.49D-12    700    15.00
 RR100   1400  700    svd lapack   0.86D-14  0.91D-12    700    84.00
 RR100   1400  700    gqr lapack   0.36D-14  0.16D-11    700    22.00
condition number:     0.49D+03

 RR100   2000  400    huang6       0.18D-11  0.40D-14    400    11.00
 RR100   2000  400    mod.huang6   0.55D-14  0.49D-15    400    22.00
 RR100   2000  400    huang7       0.18D-11  0.41D-14    400    11.00
 RR100   2000  400    mod.huang7   0.21D-14  0.12D-14    400    21.00
 RR100   2000  400    impl.qr5     0.60D-14  0.33D-15    400    18.00
 RR100   2000  400    qr lapack    0.30D-14  0.59D-14    400     8.00
 RR100   2000  400    svd lapack   0.84D-14  0.74D-14    400    22.00
 RR100   2000  400    gqr lapack   0.27D-14  0.15D-13    400    11.00
conition number:      0.26D+03

 IDF1    1050  950    huang6       0.17D-03  0.25D-14    950    32.00
 IDF1    1050  950    mod.huang6   0.44D-10  0.23D-16    950    63.00
 IDF1    1050  950    huang7       0.17D-03  0.98D-15    950    30.00
 IDF1    1050  950    mod.huang7   0.24D-10  0.23D-17    950    60.00
 IDF1    1050  950    impl.qr5     0.32D-07  0.42D-12    950    58.00
 IDF1    1050  950    qr lapack    0.23D-09  0.14D-14    950    17.00
 IDF1    1050  950    svd lapack   0.22D-09  0.16D-14    950   128.00
 IDF1    1050  950    gqr lapack   0.18D-09  0.52D-14    950    29.00
condition number:     0.15D+08
\end{verbatim}           

\newpage

Test results for overdetermined linear systems - continued

\begin{verbatim}           
 matrix  dimension    method       solution  residual   rank     time
	   m    n                    error     error
 --------------------------------------------------------------------

 IDF1    1400  700    huang6       0.53D-03  0.21D-14    700    23.00
 IDF1    1400  700    mod.huang6   0.61D-10  0.11D-14    700    47.00
 IDF1    1400  700    huang7       0.53D-03  0.11D-14    700    24.00
 IDF1    1400  700    mod.huang7   0.23D-10  0.45D-15    700    47.00
 IDF1    1400  700    impl.qr5     0.28D-07  0.30D-12    700    43.00
 IDF1    1400  700    qr lapack    0.24D-09  0.23D-14    700    15.00
 IDF1    1400  700    svd lapack   0.23D-09  0.89D-15    700    63.00
 IDF1    1400  700    gqr lapack   0.23D-09  0.33D-14    700    23.00
condition number:     0.18D+08

 IDF1    2000  400    huang6       0.12D-02  0.17D-14    400    11.00
 IDF1    2000  400    mod.huang6   0.14D-09  0.90D-15    400    22.00
 IDF1    2000  400    huang7       0.12D-02  0.21D-15    400    11.00
 IDF1    2000  400    mod.huang7   0.31D-10  0.98D-15    400    21.00
 IDF1    2000  400    impl.qr5     0.20D-07  0.18D-12    400    17.00
 IDF1    2000  400    qr lapack    0.36D-09  0.97D-15    400     8.00
 IDF1    2000  400    svd lapack   0.36D-09  0.87D-15    400    17.00
 IDF1    2000  400    gqr lapack   0.52D-09  0.82D-16    400    12.00
condition number:     0.33D+08

 IDF2    1050  950    huang6       0.65D+03  0.18D-12    950    29.00
 IDF2    1050  950    mod.huang6   0.83D+10  0.86D-05      4     1.00
 IDF2    1050  950    huang7       0.65D+03  0.76D-13    950    30.00
 IDF2    1050  950    mod.huang7   0.83D+10  0.17D-05      4     0.00
 IDF2    1050  700    impl.qr5     --- break-down ---   
 IDF2    1050  950    qr lapack    0.11D+14  0.11D-05    950    16.00
 IDF2    1050  950    svd lapack   0.10D+01  0.36D-15      3   138.00
 IDF2    1050  950    gqr lapack   0.10D+01  0.31D-15      3    25.00
 condition number:    0.11D+21
 
 IDF2    1400  700    huang6       0.18D+04  0.97D-12    700    20.00
 IDF2    1400  700    mod.huan     0.23D+11  0.50D-04      4     1.00
 IDF2    1400  700    huang7       0.18D+04  0.38D-12    700    20.00
 IDF2    1400  700    mod.huang7   0.23D+11  0.35D-04      4     0.00
 IDF2    1050  700    impl.qr5     --- break-down --- 
 IDF2    1400  700    qr lapack    0.32D+13  0.41D-06    700    14.00
 IDF2    1400  700    svd lapack   0.10D+01  0.29D-15      3    74.00
 IDF2    1400  700    gqr lapack   0.10D+01  0.68D-15      3    21.00
 condition number:    0.45D+20
\end{verbatim}           

\newpage

Test results for overdetermined linear systems - continued

\begin{verbatim}           
 matrix  dimension    method       solution  residual   rank     time
	   m    n                    error     error
 --------------------------------------------------------------------

 IDF2    2000  400    huang6       0.10D+04  0.51D-12    400     9.00
 IDF2    2000  400    mod.huang6   0.76D+10  0.33D-05      4     0.00
 IDF2    2000  400    huang7       0.10D+04  0.50D-12    400     9.00
 IDF2    2000  400    mod.huang7   0.76D+10  0.45D-05      4     0.00
 IDF2    1050  700    impl.qr5     --- break-down ---     
 IDF2    2000  400    qr lapack    0.91D+12  0.22D-07    400     7.00
 IDF2    2000  400    svd lapack   0.10D+01  0.68D-15      3    18.00
 IDF2    2000  400    gqr lapack   0.10D+01  0.29D-14      3    11.00
 condition number:    0.17D+20

 IDF3    1050  950    huang6       0.32D+04  0.80D-13    950    30.00
 IDF3    1050  950    mod.huang6   0.14D+04  0.14D-11      2     1.00
 IDF3    1050  950    huang7       0.32D+04  0.52D-13    950    31.00
 IDF3    1050  950    mod.huang7   0.14D+04  0.20D-09      2     0.00
 IDF3    1050  950    impl.qr      0.14D+04  0.99D-15      2     0.00
 IDF3    1050  950    qr lapack    0.37D+13  0.83D-02    950    17.00
 IDF3    1050  950    svd lapack   0.10D+01  0.24D-14      2   145.00
 IDF3    1050  950    gqr lapack   0.10D+01  0.22D-14      2    27.00
 condition number:    0.16D+21

 IDF3    1400  700    huang6       0.72D+03  0.24D-12    700    21.00
 IDF3    1400  700    mod.huang6   0.89D+03  0.17D-11      2     0.00
 IDF3    1400  700    huang7       0.72D+03  0.51D-13    700    21.00
 IDF3    1400  700    mod.huang7   0.89D+03  0.30D-09      2     0.00
 IDF3    1400  700    impl.qr5     0.89D+03  0.17D-11      2     0.00
 IDF3    1400  700    qr lapack    0.10D+13  0.17D-02    700    15.00
 IDF3    1400  700    svd lapack   0.10D+01  0.31D-13      2    76.00
 IDF3    1400  700    gqr lapack   0.10D+01  0.29D-13      2    24.00
 condition number:     0.27D+20

 IDF3    2000  400    huang6       0.38D+04  0.58D-12    400    10.00
 IDF3    2000  400    mod.huang6   0.44D+03  0.94D-16      2     0.00
 IDF3    2000  400    huang7       0.38D+04  0.35D-12    400     9.00
 IDF3    2000  400    mod.huang7   0.44D+03  0.67D-12      2     0.00
 IDF3    2000  400    impl.qr5     0.44D+03  0.62D-16      2     0.00
 IDF3    2000  400    qr lapack    0.45D+12  0.24D-02    400     8.00
 IDF3    2000  400    svd lapack   0.10D+01  0.65D-15      2    17.00
 IDF3    2000  400    gqr lapack   0.10D+01  0.19D-14      2    12.00
 condition number:    0.63D+19
\end{verbatim}           

\newpage

Test results for overdetermined linear systems - continued

\begin{verbatim}
 matrix  dimension    method       solution  residual   rank     time
	   m    n                    error     error
 --------------------------------------------------------------------

 IR50    1050  950    huang6       0.10D+02  0.34D-15    950    33.00
 IR50    1050  950    mod.huang6   0.92D+01  0.35D-13    773    61.00
 IR50    1050  950    huang7       0.10D+02  0.66D-14    950    32.00
 IR50    1050  950    mod.huang7   0.92D+01  0.48D-13    773    58.00
 IR50    1050  950    impl.qr5     0.92D+01  0.68D-13    773    52.00
 IR50    1050  950    qr lapack    0.44D+12  0.16D-01    950    17.00
 IR50    1050  950    svd lapack   0.41D+00  0.13D-13    773    96.00
 IR50    1050  950    gqr lapack   0.41D+00  0.17D-13    773    28.00
condition number:     0.65D+21

 IR50    1400  700    huang6       0.22D+02  0.29D-13    700    24.00
 IR50    1400  700    mod.huang6   0.64D+01  0.12D-12    618    48.00
 IR50    1400  700    huang7       0.22D+02  0.26D-13    700    24.00
 IR50    1400  700    mod.huang7   0.64D+01  0.18D-13    618    48.00
 IR50    1400  700    impl.qr5     0.64D+01  0.57D-13    618    37.00
 IR50    1400  700    qr lapack    0.34D+12  0.87D-01    700    15.00
 IR50    1400  700    svd lapack   0.36D+00  0.62D-13    618    50.00
 IR50    1400  700    gqr lapack   0.36D+00  0.11D-12    618    22.00
condition number:     0.43D+21

 IR50    2000  400    huang6       0.41D+04  0.30D-11    400    10.00
 IR50    2000  400    mod.huang6   0.45D+00  0.12D-14    374    20.00
 IR50    2000  400    huang7       0.41D+04  0.71D-12    400    11.00
 IR50    2000  400    mod.huang7   0.45D+00  0.24D-14    374    21.00
 IR50    2000  400    impl.qr5     0.45D+00  0.43D-14    374    17.00
 IR50    2000  400    qr lapack    0.80D+12  0.19D+00    400     7.00
 IR50    2000  400    svd lapack   0.33D+00  0.14D-13    374    18.00
 IR50    2000  400    gqr lapack   0.33D+00  0.11D-12    374    11.00
condition number:     0.11D+21

\end{verbatim} 

\end{document}